\documentclass[11pt]{article}
 
\usepackage{epsfig}
\usepackage{latexsym}
\usepackage{amsfonts}
\usepackage{amssymb}

\newtheorem{definition}{Definition}[section]
\newtheorem{proposition}[definition]{Proposition}
\newtheorem{lemma}[definition]{Lemma}
\newtheorem{theorem}[definition]{Theorem}
\newtheorem{corollary}[definition]{Corollary}

\begin{document}

\baselineskip=7mm

\begin{center}

\vspace*{10mm}

{\Large \bf Connectivity Preserving Hamiltonian Cycles in \\ $k$-Connected Dirac Graphs} 

\vspace*{30mm}

{\large Toru Hasunuma}

\vspace*{8mm}

{Department of Mathematical Sciences \\
Tokushima University \\
2--1 Minamijosanjima, Tokushima, 770--8506, Japan \\
}

\vspace*{5mm}
{\tt hasunuma@tokushima-u.ac.jp}

\end{center}

\vspace*{20mm}

\noindent\textbf{Abstract:}
We show that for $k \geq 2$, there exists a function $f(k) = O(k)$ such that 
every $k$-connected graph $G$ of order $n \geq f(k)$ with minimum degree at least
$\frac{n}{2}$ contains a Hamiltonian cycle $H$ such that $G-E(H)$ is $k$-connected. 
Applying Nash-Williams' result on edge-disjoint Hamiltonian cycles, 
we also show that for $k \geq 2$ and $\ell \geq 2$, 
there exists a function $g(k,\ell) = O(k\ell)$ such that 
every $k$-connected graph $G$ of order $n \geq g(k,\ell)$ 
with minimum degree at least $\frac{n}{2}$ contains $\ell$ edge-disjoint Hamiltonian cycles 
$H_1,H_2,\ldots,H_\ell$ such that 
$G-\cup_{1 \leq i \leq \ell}E(H_i)$ is $k$-connected.
As a corollary, we have a statement that refines the result of Nash-Williams 
for $k$-connected graphs with $k \leq 8$.
Moreover, when the connectivity of $G$ is exactly $k$, a similar result with an improved lower bound
on $n$ can be shown, which does not depend on the result of Nash-Williams.

\vspace*{10mm}


\newpage

\bigskip

\section{Introduction}

Throughout this paper, a graph $G = (V,E)$ means a simple undirected graph. 
For $v \in V(G)$, $N_G(v)$ denotes 
the set of vertices adjacent to $v$ in $G$. 
The degree of a vertex $v$ in $G$ is denoted by ${\rm deg}_G(v)$, i.e.,
${\rm deg}_G(v) = |N_G(v)|$.
Let $\delta(G) = \min_{v \in V(G)}{\rm deg}_G(v)$ and 
$\Delta(G) = \max_{v \in V(G)}{\rm deg}_G(v)$.
A graph is $k$-connected if the connectivity $\kappa(G)$ of $G$ is at least $k$. 
For a proper subset $S \subsetneq V(G)$ (respectively, a subset $F \subseteq E(G)$), 
we denote by $G-S$ (respectively, $G-F$) the graph obtained from $G$ 
by deleting every vertex in $S$ (respectively, every edge in $F$),
where $G-\{a\}$ is abbreviated to $G-a$.
For two sets $A$ and $B$, $A \setminus B$ denotes the set difference 
$\{ x\ |\ x \in A, x \not\in B \}$.
For a nonempty subset $S \subseteq V(G)$, the subgraph of $G$ induced by $S$ is denoted
by $\langle S \rangle_G$, i.e., $\langle S \rangle_G = G-(V(G) \setminus S)$.
For $F \subseteq E(G)$, the edge-induced subgraph of $G$ by $F$ is also denoted by 
$\langle F \rangle_G$, i.e.,
$V(\langle F \rangle_G) = \{ u \in V(G)\ |\ uv \in F\}$ and 
$E(\langle F \rangle_G) = F$.

In 1974, Mader showed the following result on the existence of a connectivity preserving cycle
in $k$-connected graphs.

\begin{theorem} {\rm (Mader \cite{M}) } \label{M}
Every $k$-connected graph $G$ with $\delta(G) \geq k+2$ contains a cycle $C$
such that $G-E(C)$ is $k$-connected.
\end{theorem}

Note that the length of a connectivity preserving cycle $C$ is not specified in Theorem \ref{M}.
When $k = 2$, Jackson \cite{J} proved a stronger result 
that every 2-connected graph $G$ with $e \in E(G)$ and $\delta(G) \geq 4$
has a cycle $C$ of length at least $\delta(G)-1$ 
such that $e \not\in E(C)$ and $G-E(C)$ is 2-connected. 
Apart from these results, connectivity preserving trees in $k$-connected or $k$-edge-connected
graphs have recently been studied in \cite{H}.

In 1952, Dirac presented the following famous sufficient condition for a graph to be Hamiltonian.

\begin{theorem} {\rm (Dirac \cite{D})} \label{D}
Every graph $G$ of order $n \geq 3$ with $\delta(G) \geq \frac{n}{2}$ is Hamiltonian.
\end{theorem}

Any graph $G$ of order $n \geq 3$ with $\delta(G) \geq \frac{n}{2}$ is often called
a Dirac graph after Theorem \ref{D}. 

In 1971, Nash-Williams showed that Theorem \ref{D} can be strengthened to
the following result on edge-disjoint Hamiltonian cycles under a lower bound condition on $n$.

\begin{theorem} {\rm (Nash-Williams \cite{N})} \label{N}
For any $\ell \geq 2$, 
every graph $G$ of order $n \geq \frac{224\ell}{5}-10$ with $\delta(G) \geq \frac{n}{2}$ contains 
$\ell$ edge-disjoint Hamiltonian cycles.
\end{theorem}

In this paper, we strengthen Theorem \ref{D} from the connectivity preserving point of view
with a lower bound condition on $n$; that is, we show the following.

\begin{theorem} \label{main1}
For any $k \geq 2$, 
every $k$-connected graph $G$ of order $n \geq 6k+12-2\left\lceil \frac{k-2}{k} \right\rceil$ 
with $\delta(G) \geq \frac{n}{2}$
contains a Hamiltonian cycle $H$ such that $G-E(H)$ is $k$-connected.
\end{theorem}

Applying Theorem \ref{N}, 
we also show the following result on 
connectivity preserving edge-disjoint Hamiltonian cycles
in Dirac graphs.

\begin{theorem} \label{main2}
For any $k \geq 2$ and any $\ell \geq 2$, 
every $k$-connected graph $G$ of order 
$n \geq \max\left\{k\ell+\max\{k\ell, 6\ell+2 \}+3k+2\ell-6,\ 6k+20\ell-8
-2\left\lceil \frac{k-2}{k} \right\rceil,\ \frac{224\ell}{5}-10\right\}$ 
with $\delta(G) \geq \frac{n}{2}$
contains $\ell$ edge-disjoint Hamiltonian cycles 
$H_1,H_2,\ldots,H_\ell$ such that 
$G-\cup_{1 \leq i \leq \ell}E(H_i)$ is $k$-connected.
\end{theorem}

When $k \leq 8$, 
the lower bound on $n$ in Theorem \ref{main2} is the same as 
that in Theorem \ref{N}. 
Thus, we have the following corollary, which refines Theorem \ref{N}. 

\begin{corollary} \label{main-cor}
For any $2 \leq k \leq 8$ and any $\ell \geq 2$, 
every $k$-connected graph $G$ of order $n \geq \frac{224\ell}{5}-10$ 
with $\delta(G) \geq \frac{n}{2}$ contains 
$\ell$ edge-disjoint Hamiltonian cycles
$H_1,H_2,\ldots,H_\ell$ such that 
$G-\cup_{1 \leq i \leq \ell}E(H_i)$ is $k$-connected.
\end{corollary}

Moreover, we show that if $\kappa(G) = k$, then 
the lower bounds on $n$ in Theorems \ref{main1} and \ref{main2} 
can be improved as follows.

\begin{theorem} \label{main3}
For any $\ell \geq 1$, every graph $G$ of order 
$n \geq \max\{2\kappa(G)\ell+\kappa(G)-2\ell+1,\ 6\kappa(G)+8\ell-4
-2\left\lceil \frac{\kappa(G)-2}{\kappa(G)} \right\rceil \}$
with $\delta(G) \geq \frac{n}{2}$
contains $\ell$ edge-disjoint Hamiltonian cycles 
$H_1,H_2,\ldots,H_\ell$ such that 
$\kappa(G-\cup_{1 \leq i \leq \ell}E(H_i)) = \kappa(G)$. 
\end{theorem}

Note that unlike Theorem \ref{main2}, Theorem \ref{main3} does not depend on Theorem \ref{N}.  
In particular, when $\kappa(G) \leq 11$,
the lower bound on $n$ in Theorem \ref{main3} improves that in Theorem \ref{N}.

Section 2 presents results on edge-disjoint pairs of independent edges.
Sufficient conditions for a graph to have edge-disjoint
Hamiltonian paths connecting any specified pairs of vertices are given in Section 3. 
Section 4 presents a proof of a more general statement than Theorem \ref{main1},
from which Theorem \ref{main2} is also obtained by applying Theorem \ref{N}. 
Modifying the proof of the statement in Section 4, Theorem \ref{main3} can be shown,
whose proof is given in Section 5.
Section 6 finally concludes the paper with several remarks.

\section{Edge-Disjoint Edge-Pairs}

Edges of a graph are {\it independent} if they are mutually incident with no common vertex.
In particular, we call two independent edges an {\it edge-pair}. 
A graph $G$ is decomposed into edge-pairs if $E(G)$ is partitioned into edge-pairs, i.e.,
$E(G) = E_1 \cup E_2 \cup \cdots \cup E_t$ where $E_i \cap E_j = \emptyset$ for $i \neq j$
such that each $E_i$ consists of one edge-pair.
We denote by $K_n$ the complete graph of order $n$.
The union of two vertex-disjoint graphs $G_1$ and $G_2$ is denoted by $G_1 \cup G_2$. 
We also denote by $kH$ the $k$ vertex-disjoint copies of a graph $H$. 

It is stated in \cite{R} that the following (unpublished) result is due to Caro. 
It is also mentioned in \cite{CL} that the proof is quite lengthy. 
For the sake of readability (and also for constructions of Hamiltonian cycles in our main results), 
we present our original proof for the result, which might be simpler than that by Caro.

\begin{theorem} \label{decom}
A nonempty graph $G$ can be decomposed into edge-pairs if and only if 
the following three conditions hold: 
\begin{enumerate}
\item $|E(G)|$ is even,
\item $|E(G)| \geq 2 \Delta(G)$,
\item $K_3 \not\subset G$ when $|E(G)| = 4$.
\end{enumerate}
\end{theorem}

\noindent {\bf Proof:} 
Suppose that $G$ is decomposed into edge-pairs $E_1, E_2, \ldots, E_t$.
Then $|E(G)|$ is clearly even.
Let $x$ be a vertex of degree $\Delta(G)$.
For any edge $e$ incident with $x$, another edge $e'$ 
in the edge-pair $E_i$ containing $e$ is not incident with $x$.
Thus, $t \geq \Delta(G)$ and $G$ has at least $2\Delta(G)$ edges.
Since no two edges in a triangle are independent, $G$ has no triangle 
when $|E(G)| = 4$. 

Conversely, suppose that all the three conditions hold.
When $\Delta(G) = 1$, $G$ is trivially decomposed into edge-pairs since $|E(G)|$ is even.
Suppose that $\Delta(G) \geq 2$.
We proceed by induction on $\Delta(G)$ to show that $G$ is decomposed into edge-pairs
when $|E(G)| = 2\Delta(G)$.

\medskip

\noindent Case 1: $G$ has exactly one vertex of degree $\Delta(G)$.

Let $x$ be the vertex of degree $\Delta(G)$.
Since $G \not\cong K_3 \cup (n-3)K_1$ where $n \geq 3$ 
and $G \not\cong S_p \cup (n-p)K_1$ where $S_p$ is a star of order $p$ where $n \geq p$,
we can select an edge-pair $\{xy,x'y'\}$ in $G$. 
In particular, when $|E(G)| = 6$, i.e., $\Delta(G) = 3$, 
$G$ has at most one triangle and thus we can select an edge-pair $\{xy,x'y'\}$  
so that $G-\{xy,x'y'\}$ has no triangle. 
Since $|E(G-\{xy,x'y'\})| = 2\Delta(G-\{xy,x'y'\})$, 
by the induction hypothesis, we have a decomposition of $G$ into edge-pairs.

\medskip

\noindent Case 2: $G$ has exactly two vertices of degree $\Delta(G)$.

Let $x_1,x_2$ be the vertices of degree $\Delta(G)$.
If $x_1x_2 \not\in E(G)$ (respectively, $x_1x_2 \in E(G)$), 
then we can select an edge-pair $\{x_1x'_1, x_2x'_2\}$ (respectively, $\{x_1x_2, x'_1x'_2\}$)
such that $|E(G-\{x_1x'_1,x_2x'_2\})| = 2\Delta(G-\{x_1x'_1,x_2x'_2\})$ 
(respectively, $|E(G-\{x_1x_2,x'_1x'_2\})| = 2\Delta(G-\{x_1x_2,x'_1x'_2\})$).
Note that if $x_1x_2 \not\in E(G)$ (respectively, $x_1x_2 \in E(G)$), 
then $G$ (respectively, $G-\{x_1x_2,x'_1x'_2\}$) is bipartite and has no triangle.
Thus, we can apply the induction hypothesis to obtain a decomposition of $G$ into edge-pairs.

\medskip

\noindent Case 3: $G$ has at least three vertices of degree $\Delta(G)$.

Let $x_1,x_2,x_3$ be vertices of degree $\Delta(G)$.
Since $|E(G)| \geq 3\Delta(G)-|E(\langle \{x_1,x_2,x_3\} \rangle_G)|$, 
it holds that $\Delta(G) \leq |E(\langle \{x_1,x_2,x_3\} \rangle_G)|$.
By the third condition, it does not happen that
$\Delta(G) = 2$ and $|E(\langle \{x_1,x_2,x_3\} \rangle_G)| = 3$.
If $\Delta(G) = |E(\langle \{x_1,x_2,x_3\} \rangle_G)| = 2$, then 
$G$ is either a path of order 5 or a cycle of order 4, which implies that
$G$ can be decomposed into two edge-pairs.
If $\Delta(G) = |E(\langle \{x_1,x_2,x_3\} \rangle_G)| = 3$,
then $G$ has the triangle $\langle \{x_1,x_2,x_3\} \rangle_G$ such that
each vertex $x_i$ is adjacent to exactly one vertex $x'_i$ not in the triangle.
Thus, in such a case, $G$ is decomposed into three edge-pairs
$\{x_1x'_1,x_2x_3\}, \{x_2x'_2,x_1x_3\}, \{x_3x'_3,x_1x_2\}$.

\medskip

Let $\Phi(G) = |E(G)|-2\Delta(G)$.
By the previous argument, we have shown the case that $\Phi(G) = 0$.
Now we employ double induction on
$\Delta(G)$ and $\Phi(G)$. 
Suppose that $\Delta(G) \geq 2$ and $\Phi(G) \geq 2$.
Let $\{uv,xy\}$ be an edge-pair in $G$.
Note that when $|E(G)| = 6$, any triangle in $G$ (if exists) forms
a connected component since $\Delta(G) = 2$.
Thus, in such a case, we can select an edge-pair $\{uv,xy\}$ so that $G-\{uv,xy\}$ has no triangle. 
Let $H = G - \{uv,xy\}$.
Then, either $\Delta(H) = \Delta(G)-1$ or $\Delta(H) = \Delta(G)$.
If $\Delta(H) = \Delta(G)-1$, then $\Phi(H) = \Phi(G)$.
If $\Delta(H) = \Delta(G)$, then $\Phi(H) = \Phi(G)-2$.
Therefore, in either case, we can apply the double induction hypothesis.
Hence, we have the desired result.
\hfill $\blacksquare$

\bigskip

Based on Theorem \ref{decom}, the maximum number of edge-disjoint edge-pairs in a graph 
is determined as follows.

\begin{theorem} \label{max-DEP}
For any graph $G$ of order $n$ such that $G \not\cong K_3 \cup (n-3)K_1$
for $n \geq 3$
and $G \not\cong K_3 \cup K_2 \cup (n-5)K_1$ for $n \geq 5$, 
the maximum number of edge-disjoint edge-pairs in $G$ is given by 
$$\displaystyle \left\{\begin{array}{ll}
\left\lfloor \frac{|E(G)|}{2} \right\rfloor & \mbox{ if } |E(G)| \geq 2\Delta(G), \\[4mm]
|E(G)|-\Delta(G) & \mbox{ if }  |E(G)| < 2\Delta(G).
\end{array} \right.$$
\end{theorem}

\noindent {\bf Proof:}
Let $S$ be a maximum set of edge-disjoint edge-pairs in $G$.
Clearly, $|S| \leq \left\lfloor \frac{|E(G)|}{2} \right\rfloor$. 
Let $x$ be a vertex of degree $\Delta(G)$ in $G$. 
For any edge-pair in $S$, it must contain an edge not incident with $x$.
Thus, $|S| \leq |E(G)|-\Delta(G)$. 

When $\Delta(G) \leq 1$, the statement clearly holds.
Suppose that $\Delta(G) \geq 2$.

\medskip

\noindent Case 1: $|E(G)| \geq 2\Delta(G)$. 

If $|E(G)|$ is even, then by Theorem \ref{decom}, 
$G$ can be decomposed into $\frac{|E(G)|}{2}$ disjoint edge-pairs since 
$G \not\cong K_3 \cup K_2 \cup (n-5)K_1$ for $n \geq 5$.
Suppose that $|E(G)|$ is odd.
Then for any $e \in E(G)$, $|E(G-e)| \geq 2\Delta(G-e)$. 
In particular, 
when $|E(G)| = 5$, i.e., $\Delta(G) = 2$, we can select an edge $e$ so that $K_3 \not\subset G-e$. 
Therefore, by Theorem \ref{decom}, $G-e$ can be decomposed into $\frac{|E(G)|-1}{2}$
edge-pairs.

\medskip

\noindent Case 2: $|E(G)| < 2\Delta(G)$.

If $\Delta(G) = |E(G)| = 2$, then $G$ has no edge-pair.
If $\Delta(G) = 2$ and $|E(G)| = 3$, then 
$G$ has one edge-pair since $G \not\cong K_3 \cup (n-3)K_1$ for $n \geq 3$. 
Moreover, if $G$ has at least three vertices of degree $\Delta(G)$, then 
$\Delta(G) = 2$ and $|E(G)| = 3$ and thus $G$ has an edge-pair. 

Suppose that $\Delta(G) \geq 3$ and 
$G$ has at most two vertices of degree $\Delta(G)$. 

\medskip

\noindent Case 2.1: $G$ has exactly two vertices of degree $\Delta(G)$.

Let $x,y$ be the vertices of degree $\Delta(G)$. 
Then $x$ and $y$ must be adjacent and $|E(G)| = 2\Delta(G)-1$.
Note that $G-xy$ is bipartite and has no triangle. 
Since $|E(G-xy)| = 2\Delta(G-xy)$, by Theorem \ref{decom}, 
there are $\Delta(G)-1 = |E(G)| - \Delta(G)$ edge-disjoint edge-pairs in $G-xy \subset G$. 

\medskip

\noindent Case 2.2: $G$ has exactly one vertex of degree $\Delta(G)$.

Let $x$ be the vertex of degree $\Delta(G)$.
Let $G'$ be a graph obtained from $G$ by deleting
$2\Delta(G)-|E(G)|$ edges incident with $x$ in such a way that
if $\Delta(G-x) = |E(G)|-\Delta(G)$ and $G$ has an edge joining
$x$ and a vertex $z$ of degree $\Delta(G-x)$ in $G-x$, then the edge $xz$ is preferentially deleted.
We consequently have ${\rm deg}_{G'}(x) = \Delta(G') = |E(G)|-\Delta(G)$ and 
$|E(G')| = 2\Delta(G')$.
If $G' \cong K_3 \cup K_2 \cup (n-5)K_1$ for $n \geq 5$, 
then the triangle in $G'$ has the vertex $x$ and 
we can modify $G'$ by replacing 
an edge incident with $x$ in the triangle with an edge in $E(G) \setminus E(G')$.
Therefore, by Theorem \ref{decom}, there are $(|E(G)|-\Delta(G))$ edge-disjoint edge-pairs 
in $G' \subset G$.
\hfill $\blacksquare$

\bigskip

\section{Edge-Disjoint Hamiltonian Paths}

In 1962, P\'{o}sa presented the following sufficient condition for a graph to be Hamiltonian,
which extends Dirac's condition. 

\begin{theorem} {\rm (P\'{o}sa \cite{P})} \label{P}
Let $G$ be a graph of order $n \geq 3$ such that for every integer $j$ with 
$1 \leq j < \frac{n}{2}$, the number of vertices of degree not exceeding $j$ is less than $j$.
Then $G$ is Hamiltonian. 
\end{theorem}

A graph $G$ is {\it Hamiltonian-connected} if for every pair of vertices in $G$,
there exists a Hamiltonian path connecting them.

In 1969, Chartrand, Kapoor and Kronk showed 
the following sufficient condition for a graph to be 
Hamiltonian-connected, which is analogous to Theorem \ref{P}.

\begin{theorem} {\rm (Chartrand, Kapoor, Kronk \cite{CKK})} \label{CKK}
Let $G$ be a graph of order $n \geq 4$ such that for every integer $j$ with 
$2 \leq j \leq \frac{n}{2}$, the number of vertices of degree not exceeding $j$ is less than $j-1$.
Then $G$ is Hamiltonian-connected.
\end{theorem}

For a graph $G$ and a nonnegative integer $j$, 
let 
$$\psi_j(G) = |\{ v \in V(G)\ |\ {\rm deg}_G(v) \leq j \}|.$$

Theorem \ref{CKK} can be naturally extended to edge-disjoint Hamiltonian paths connecting
any specified pairs of vertices as follows.
Note that the same pair of vertices may appear twice or more in a multiset of $\ell$ pairs of vertices. 

\begin{proposition} \label{EDHP}
Let $G$ be a graph of order $n \geq 4\ell$ such that 
for every integer $j$ with 
$2\ell \leq j \leq \frac{n}{2}+2(\ell-1)$, 
$\psi_j(G) < j-2\ell+1$.
Then for any multiset 
$\{ \{u_1,v_1\}, \{u_2,v_2\}, \ldots, \{u_\ell,v_\ell \} \}$
of $\ell$ pairs of vertices in $G$, 
there are $\ell$ edge-disjoint Hamiltonian paths $P_1,P_2,\ldots,P_\ell$ in $G$ 
such that $P_i$ connects $u_i$ and $v_i$ for each $1 \leq i \leq \ell$.
\end{proposition}

\noindent {\bf Proof:} 
We proceed by induction on $\ell$.
When $\ell = 1$, the statement is the same as Theorem \ref{CKK}.
Suppose that $\ell \geq 2$. 
Let $G$ be a graph of order $n \geq 4\ell$ such that 
for every $j$ with 
$2\ell \leq j \leq \frac{n}{2}+2(\ell-1)$, 
$\psi_j(G) < j-2\ell+1$.
Let $\{ \{u_1,v_1\}, \{u_2,v_2\}, \ldots, \{u_\ell,v_\ell \} \}$ be 
a multiset of $\ell$ pairs of vertices in $G$. 

Since $\psi_j(G) = 0$ for any $j \leq 2\ell$ and 
$\psi_j(G) < j-3$ for any $2\ell < j \leq \frac{n}{2}$, 
by Theorem \ref{CKK}, there exists a Hamiltonian path $P$ connecting $u_{\ell}$ and $v_{\ell}$ in $G$. 
Let $G' = G-E(P)$.
Then it holds that ${\rm deg}_{G'}(u_{\ell}) = {\rm deg}_{G}(u_{\ell}) -1$,
${\rm deg}_{G'}(v_{\ell}) = {\rm deg}_{G}(v_{\ell}) -1$,
and for any vertex $v \in V(G') \setminus \{u_{\ell},v_{\ell}\}$, ${\rm deg}_{G'}(v) = {\rm deg}_G(v)-2$.
This means that for each $j \geq 2$, 
$$\psi_{j-2}(G') \leq \psi_{j}(G).$$
Thus, for every $j$ with $2\ell \leq j \leq \frac{n}{2}+2(\ell-1)$,
$$\psi_{j-2}(G') < j-2\ell+1,$$ 
that is, for every $j$ with $2(\ell-1) \leq j \leq \frac{n}{2}+2(\ell-2)$,
$$\psi_{j}(G') < j-2(\ell-1)+1.$$ 
Therefore, by the induction hypothesis, 
there are $\ell-1$ edge-disjoint Hamiltonian paths $P_1,P_2,\ldots,P_{\ell-1}$ in $G'$ 
such that $P_i$ connects $u_i$ and $v_i$ for each $1 \leq i < \ell$, 
and thus the proof is complete.
\hfill $\blacksquare$

\bigskip

From Proposition \ref{EDHP}, the following lemma is obtained.

\begin{lemma} \label{H}
Let $G$ be a graph of order $n$.
Suppose that $V(G)$ is partitioned into $V_1$ and $V_2$ such that
every vertex in $V_1$ has degree at least $\frac{n}{2}+2\ell-1$ and 
every vertex in $V_2$ has degree at least $|V_2|+2\ell$, 
where $V_1$ is nonempty and $V_2$ is possibly empty.
Then for any multiset 
$\{ \{u_1,v_1\}, \{u_2,v_2\}, \ldots, \{u_\ell,v_\ell \} \}$
of $\ell$ pairs of vertices in $G$, 
there are $\ell$ edge-disjoint Hamiltonian paths $P_1,P_2,\ldots,P_\ell$ in $G$ 
such that $P_i$ connects $u_i$ and $v_i$ for each $1 \leq i \leq \ell$.
\end{lemma}

\noindent {\bf Proof:} 
Since $V_1 \neq \emptyset$, $\Delta(G) \geq \frac{n}{2}+2\ell-1$ and we have $n \geq 4\ell$.
If $|V_2| \geq \frac{n}{2}-1$, then $\delta(G) \geq \frac{n}{2}+2\ell-1$ and 
$G$ trivially satisfies the condition in Proposition \ref{EDHP}. 
Suppose that $|V_2| \leq \frac{n}{2}-2$. 
Since every vertex in $V_1$ has degree at least $\frac{n}{2}+2\ell-1$, 
for every $j \leq \frac{n}{2}+2(\ell-1)$,
$$\psi_j(G) \leq |V_2|.$$
Since every vertex in $V_2$ has degree at least $|V_2|+2\ell$,
for every $j < |V_2|+2\ell \leq \frac{n}{2}+2(\ell-1)$,
$$\psi_j(G) = 0.$$
On the other hand, for every $j \geq |V_2|+2\ell$,
$$|V_2| < j-2\ell+1.$$
Therefore, for every $j$ with $2\ell \leq j \leq \frac{n}{2}+2(\ell-1)$,
$$\psi_j(G) < j-2\ell+1.$$
Hence, the lemma follows from Proposition \ref{EDHP}. 
\hfill $\blacksquare$

\bigskip

\section{Proofs of Theorems \ref{main1} and \ref{main2}}

In 1968, Chartrand and Harary originally presented the following fundamental sufficient condition for a graph 
to be $k$-connected, which is the same as Dirac's condition when $k = 2$. 

\begin{theorem} {\rm (Chartrand and Harary \cite{CH})} \label{CH}
For any $k \geq 1$, every graph of order $n \geq k+1$ with 
$\delta(G) \geq \frac{n+k-2}{2}$ is $k$-connected.
\end{theorem}

The graph obtained from $2K_{\frac{n+k-1}{2}}$ where $n+k$ is odd, 
by identifying $k-1$ vertices has minimum degree $\frac{n+k-3}{2}$ 
and is not $k$-connected.
In this sense, the lower bound of $\frac{n+k-2}{2}$ on $\delta(G)$ in Theorem \ref{CH} is best possible.

A vertex of degree 2 in a path $P$ is called an internal vertex in $P$ and
the set of internal vertices in $P$ is denoted by $V_I(P)$.
The following lemma 
gives a construction of a $k$-connected graph from two vertex-disjoint $k$-connected graphs.

\begin{lemma} \label{k-connect}
Let $G_1$ and $G_2$ be $k$-connected graphs such that $V(G_1) \cap V(G_2) = \emptyset$.
Let $G$ be a graph obtained from $G_1$ and $G_2$ by adding
$k$ vertex-disjoint paths $P_1,P_2,\ldots,P_{k}$ connecting vertices in $V(G_1)$
and vertices in $V(G_2)$, where $(V(G_1) \cup V(G_2)) \cap (\cup_{1 \leq i \leq k}V_I(P_i))
= \emptyset$.
Let $G'$ be a graph obtained from $G$ by adding 
edges so that
for each $1 \leq i \leq k$, every internal vertex in $P_i$ has $k$ neighbors in 
$V(G_1) \cup V(G_2)$.
Then $G'$ is $k$-connected. 
\end{lemma}

\noindent {\bf Proof:}
Let $S \subset V(G')$ such that $|S| = k-1$.
Let $S_1 = S \cap V(G_1)$ and $S_2 = S \cap V(G_2)$. 
Since $G_1$ (respectively, $G_2$) is $k$-connected, 
$G_1 - S_1$ (respectively, $G_2-S_2$) is connected.
For any $u \in V(G_1) \setminus S_1$ and any $v \in V(G_2) \setminus S_2$,
there is a path connecting them in $G'-S$ since at least one path $P_i$ 
still exists in $G'-S$.
Moreover, in $G'-S$, 
any $w \in (\cup_{1 \leq i \leq k}V_I(P_i)) \setminus S$ 
has at least one neighbor in $(V(G_1) \cup V(G_2)) \setminus S$.
Therefore, $G'-S$ is connected. 
Hence, $G'$ is $k$-connected.
\hfill $\blacksquare$

\bigskip

In what follows, we show the following.

\begin{theorem} \label{main0}
For any $k \geq 2$, $p \geq 1$ and $q \geq 1$,
every $k$-connected graph $G$ of order 
$n \geq \max\{6(k+p)+8(q-1)-2\left\lceil \frac{k-2}{k} \right\rceil, \ 
k(q+3)+\max\{(k-2)q-2,2p\}+2\max\{p,q\}-4\}$
with 
$\delta(G) \geq \frac{n}{2}$ satisfies one of the following properties:
\begin{itemize}
\item For any spanning subgraph $H$ of $G$ with $\Delta(H) \leq p$, 
$G-E(H)$ is $k$-connected.
\item There are $q$ edge-disjoint Hamiltonian cycles
$H_1,H_2,\ldots,H_q$ in $G$ such that
$G-\cup_{1\leq i \leq q}E(H_i)$ is $k$-connected.
\end{itemize}
\end{theorem}

\noindent {\bf Proof:}
Let $G$ be a $k$-connected graph of order $n$ with $\delta(G) \geq \frac{n}{2}$ such that 
\begin{equation} 
n \geq 6k+6p+8q-8-2\left\lceil \frac{k-2}{k} \right\rceil \label{lower-n1}
\end{equation}
and
\begin{equation} 
n \geq kq+3k+\max\{(k-2)q-2,2p\}+2\max\{p,q\}-4, \label{lower-n2}
\end{equation}
where $k,p$ and $q$ are positive integers such that $k \geq 2$.
Suppose that there exists a spanning subgraph $H$ of $G$ with $\Delta(H) \leq p$ such that 
$G-E(H)$ is not $k$-connected.
In order to prove Theorem \ref{main0}, 
it is sufficient to show that there are $q$ edge-disjoint Hamiltonian cycles
$H_1,H_2,\ldots,H_q$ in $G$ such that $G-\cup_{1 \leq i \leq q}E(H_i)$ is $k$-connected.

Let  $G' = G-E(H)$.
Also, let $k' = \kappa(G') < k$ and 
$$W = \{w_1,w_2,\ldots,w_{k'}\} \subset V(G')$$ such that $G'-W$ is disconnected.
Since $\Delta(H) \leq p$, we have
$\delta(G') \geq \frac{n}{2}-p$.

\bigskip

\noindent {\bf Claim 1.}\ 
The graph $G'-W$ has exactly two connected components $G_1$ and $G_2$
with $|V(G_1)| \leq |V(G_2)|$ such that 

\begin{equation}
\frac{n}{2}-k'-p+1 \leq |V(G_1)| \leq \frac{n-k'}{2}, \label{G1}
\end{equation}
\begin{equation}
\frac{n-k'}{2} \leq |V(G_2)| \leq \frac{n}{2}+p-1. \label{G2}
\end{equation}

\bigskip

\noindent {\bf Proof of Claim 1:}
Let $G_1,G_2,\ldots,G_t$ be the connected components of $G'-W$ such that
$|V(G_1)| \leq |V(G_2)| \leq \cdots \leq |V(G_t)|$. 
In the graph $G'$, every vertex in $V(G_i)$ for each $1 \leq i \leq t$ 
has no neighbor in $V(G) \setminus (V(G_i) \cup W)$.
Since $\delta(G') \geq \frac{n}{2}-p$, we obtain for each $1 \leq i \leq t$, 
$$|V(G_i) \cup W| \geq \frac{n}{2}-p+1. $$
Therefore, 
$$|V(G_1)| \geq \frac{n}{2}-k'-p+1.$$
Assume that $t \geq 3$.
Then we have $|V(G_1)| \leq \frac{n-k'}{3}$.
Thus, 
$$
\frac{n}{2}-k'-p+1 \leq \frac{n-k'}{3},$$
that is,
\begin{equation}
n \leq 4k'+6p-6, \label{Claim1'-1}
\end{equation}
which contradicts (\ref{lower-n1}).
Hence, it holds that $t = 2$ and 
$$|V(G_2)| \geq \frac{n-k'}{2}.$$
The upper bounds in (\ref{G1}) and (\ref{G2}) follow from the lower bounds
in (\ref{G2}) and (\ref{G1}), respectively.
$\Box$

\bigskip

\noindent {\bf Claim 2.}\   
There are $k$ vertex-disjoint paths $P^\ast_1,P^\ast_2,\ldots,P^\ast_k$ in $G$ 
connecting vertices in $V(G_1)$ and vertices in $V(G_2)$ such that
$|V(P^\ast_i)| = 2$ for each $1 \leq i \leq k-k''$, 
$3 \leq |V(P^\ast_i)| \leq 4$ for each $k-k''+1 \leq i \leq k$ and 
$\cup_{k-k''+1 \leq i \leq k}V_I(P^\ast_i) = W$,
where $\lceil \frac{k'}{2} \rceil \leq k'' \leq k'$.

\bigskip

\noindent {\bf Proof of Claim 2:}
Since $G'$ is $k'$-connected, 
there are $k'$ vertex-disjoint paths in $G'$ connecting vertices in $V(G_1)$ and vertices in $V(G_2)$. 
This means that
there are two subsets $\{x_1,x_2,\ldots,x_{k'}\} \subset V(G_1)$ and 
$\{y_1,y_2,\ldots,y_{k'}\} \subset V(G_2)$ 
such that $\{x_iw_i, w_iy_i\} \subset E(G)$ for each $1 \leq i \leq k'$.

Since $G$ is $k$-connected, there are $k$ vertex-disjoint paths $P_1,P_2,\ldots,P_k$ in $G$ connecting
vertices in $V(G_1)$ and vertices in $V(G_2)$. 
If $P_i$ has no vertex in $W$, then $P_i$ has has a subpath $P'_i = (u_i,v_i)$ of order 2,
where $u_i \in V(G_1)$ and $v_i \in V(G_2)$. 
If $P_i$ has a vertex in $W$, then $P_i$ has either a subpath $P'_i = (u_i,v_i)$ 
or $P'_i = (u_i,w_{i,1},w_{i,2},\ldots,w_{i,r_i},v_i)$, where 
$u_i \in V(G_1)$, $\{ w_{i,1},w_{i,2},\ldots,w_{i,r_i} \} \subseteq W$ and $v_i \in V(G_2)$.
Thus, without loss of generality, we may assume that
for each $1 \leq i \leq k-t$,
$|V(P'_i)| = 2$ and for each $k-t+1 \leq i \leq k$,
$|V(P'_i)| \geq 3$.
Let $\Pi' = \{P'_1,P'_2,\ldots,P'_k\}$. 
Also let $U = \{u_1,u_2,\ldots,u_k\}$ and $V = \{v_1,v_2,\ldots,v_k\}$.

In the graph $G$, any $w \in W$ has at least $\left\lceil \frac{\lceil \frac{n}{2} \rceil -k'+1}{2} \right\rceil$ 
neighbors in either $V(G_1)$ or $V(G_2)$.
From (\ref{lower-n1}), we have
$\left\lceil \frac{\lceil \frac{n}{2} \rceil -k'+1}{2} \right\rceil \geq \left\lceil \frac{3k+3-k'}{2} \right\rceil
\geq k+2$. 
Suppose that $|V(P'_i)| \geq 4$ where $k-t+1 \leq i \leq k$. 
If $w_{i,2}$ has $k+1$ neighbors in $V(G_1)$ (respectively, $V(G_2)$)
then we can select a vertex $u'_i \in (N_G(w_{i,2}) \cap V(G_1)) \setminus U$  
(respectively, $v'_i \in (N_G(w_{i,2}) \cap V(G_2)) \setminus V$) 
and obtain a path $P''_i = (u'_i,w_{i,2},\ldots,w_{i,r_i},v_i)$ (respectively, $P''_i = (u_i,w_{i,1},w_{i,2},v'_i)$) 
which is vertex-disjoint with any path in $\Pi' \setminus \{P'_i\}$.
Suppose that $|V(P'_i)| = 4$ where $k-t+1 \leq i \leq k$. 
If $w_{i,1}$ has $k+1$ neighbors in $V(G_2)$,
then we have a path $P''_i = (u_i,w_{i,1},v'_i)$, where $v'_i \in (N_G(w_{i,1}) \cap V(G_2)) \setminus V$.
From these observations, we may also assume without loss of generality that $3 \leq |V(P'_i)| \leq 4$
for any $k-t+1 \leq i \leq k$ and if $|V(P'_i)| = 4$ then $w_{i,1}$ (respectively, $w_{i,2}$) 
has at most $k$ neighbors in $V(G_2)$ (respectively, $V(G_1)$). 

Suppose that $W \setminus (\cup_{k-t+1 \leq i \leq k}V_I(P'_i)) \neq \emptyset$
and $w_j \in W \setminus (\cup_{k-t+1 \leq i \leq k}V_I(P'_i))$.

\medskip 

\noindent Case 1: $|N_G(w_j) \cap V(G_1)| \geq k+1$. 

Let $P_{w_j} = (u_{w_j},w_j,y_j)$ where $u_{w_j} \in (N_G(w_j) \cap V(G_1)) \setminus U$.
If $y_j \not\in V$, then we replace a path $P'_i$ of length 1 in $\Pi'$ with the path $P_{w_j}$. 
Note that such a path of length 1 exists in $\Pi'$ since $k' < k$. 
Suppose that there exists a path $P'_i$ with $v_i = y_j$ in $\Pi'$.
If $|V(P'_i)| = 2$, or $3 \leq |V(P'_i)| \leq 4$ and $|N_G(w_{i,r_i}) \cap V(G_2)| \leq k$, then 
we replace $P'_i$ with $P_{w_j}$.
If $3 \leq |V(P'_i)| \leq 4$ and $|N_G(w_{i,r_i}) \cap V(G_2)| \geq k+1$, then 
we modify $P'_i$ by changing the last vertex $v_i$ to $v'_i$ where 
$v'_i \in (N_G(w_{i,r_i}) \cap V(G_2)) \setminus V$
and replace a path $P'_{i'}$ of length 1 in $\Pi'$ with $P_{w_j}$.

\medskip

\noindent Case 2: $|N_G(w_j) \cap V(G_1)| \leq k$.  

Since $|N_G(w_j) \cap V(G_2)| \geq k+1$, 
let $P_{w_j} = (x_j,w_j,v_{w_j})$ where $v_{w_j} \in (N_G(w_j) \cap V(G_2)) \setminus V$.
If $x_j \not\in U$, then we replace a path $P'_i$ of length 1 in $\Pi'$ with $P_{w_j}$.
Suppose that there exists a path $P'_i$ with $u_i = x_j$ in $\Pi'$. 
If $|V(P'_i)| = 2$, or $3 \leq |V(P'_i)| \leq 4$ and $|N_G(w_{i,1}) \cap V(G_1)| \leq k$, 
we replace $P'_i$ with $P_{w_j}$;
otherwise we modify $P'_i$ by changing the first vertex $u_i$ to $u'_i$ where 
$u'_i \in (N_G(w_{i,1}) \cap V(G_1)) \setminus U$ 
and replace a path $P'_{i'}$ of length 1 in $\Pi'$ with $P_{w_j}$.

\medskip

In either case, we have a modified set $\Pi'' = \{P''_1,P''_2,\ldots,P''_k\}$ of $k$ vertex-disjoint paths 
which contains the path $P_{w_j}$.
If $W \setminus (\cup_{1 \leq i \leq k}V_I(P''_i)) \neq \emptyset$, then 
we can select a vertex in $W \setminus (\cup_{1 \leq i \leq k}V_I(P''_i))$ and 
apply a similar modification process to $\Pi''$. 
Note that the path $P_{w_j}$ in $\Pi''$ can be modified but not deleted in the process.
Therefore, iteratively applying the similar process,
we can finally obtain a desired set of $k$ vertex-disjoint paths $P^\ast_1, P^\ast_2, \ldots, P^\ast_k$; 
that is, every vertex in $W$ is used as an internal vertex of a path $P^\ast_i$ in the set
where $k-k''+1 \leq i \leq k$. 
Since each path in the set has at most two vertices in $W$,
the number $k''$ of paths containing a vertex in $W$ in the set is at least $\lceil \frac{k'}{2} \rceil$.
Hence, Claim 2 holds.
$\Box$

\bigskip

Define the subgraph $B$ of $G$ as $B = \langle  \cup_{1 \leq i \leq k}E(P^\ast_i) \rangle_G$.
Let $M_B =  \cup_{1 \leq i \leq k-k''}E(P^\ast_i)$ and  
$M_H = \{uv \in E(H)\ |\ u \in V(G_1), v \in V(G_2) \}$.
Also, let $U_1 = \{ v \in V(G_1)\ |\ N_{G}(v) \cap V(G_2) \neq \emptyset \}$
and $U_2 = \{ v \in V(G_2)\ |\ N_{G}(v) \cap V(G_1) \neq \emptyset \}$.  
Note that $M_B \subseteq M_H$ and $\max\{|U_1|,|U_2|\} \leq |M_H|$.
Moreover, let
$$q_1 = \min\{q,\max\{0,|M_H|-|M_B|-\max\{p,q\}\}\}$$ 
and $q_2 = q - q_1$.  
Note that $0 \leq q_1 \leq q$ and $0 \leq q_2 \leq q$.

\bigskip

\noindent {\bf Claim 3.}\  
If $q_1 > 0$, then there are $q_1$ edge-disjoint edge-pairs
$$\{u_1v_1, u_2v_2\}, \{u_3v_3, u_4v_4\}, \ldots, \{u_{2q_1-1}v_{2q_1-1},u_{2q_1}v_{2q_1} \} \subset
E(G) \setminus E(B),$$
where $u_i \in V(G_1)$ and $v_i \in V(G_2)$ for $1 \leq i \leq 2q_1$.

\bigskip

\noindent {\bf Proof of Claim 3:}
Suppose that $q_1 > 0$.
Then we have $q_1 = \min\{q,|M_H|-|M_B|-\max\{p,q\}\}$ which means that
$|M_H| \geq q_1+|M_B|+ \max\{p,q\}$.
Thus, $|M_H \setminus M_B| \geq q_1+ \max\{p,q\} \geq \max\{p+q_1,2q_1\}$.
Therefore, 
$\left\lfloor \frac{|M_H \setminus M_B|}{2} \right\rfloor \geq q_1$ and
$|M_H \setminus M_B|-\Delta(\langle M_H \setminus M_B \rangle_H) \geq
|M_H \setminus M_B|-\Delta(H) \geq q_1$.
Since $\langle M_H \setminus M_B \rangle_H$ is bipartite,
there is no triangle in $\langle M_H \setminus M_B \rangle_H$.
Hence, by Theorem \ref{max-DEP},
there are at least $q_1$ edge-disjoint edge-pairs in $\langle M_H \setminus M_B \rangle_H
\subset G-E(B)$. 
$\Box$

\bigskip

Define $W_1$ as follows: 
$$W_1 = \left\{ w \in W\ |\ |N_G(w) \cap V(G_1)| \geq |N_G(w) \cap V(G_2)| \right\}.$$
Let $W_2 = W \setminus W_1$.

\bigskip

\noindent {\bf Claim 4.}\   
If $q_2 > 0$, then there are 
$q_2$ edge-disjoint edge-pairs 
$$\{u_{2q_1+1}w^\ast_1,u_{2q_1+2}w^\ast_2\}, \{u_{2q_1+3}w^\ast_1,u_{2q_1+4}w^\ast_2\}, \ldots, 
\{u_{2q-1}w^\ast_1,u_{2q}w^\ast_2\} \subset E(G) \setminus E(B)$$ 
and $q_2$ edge-disjoint edge-pairs 
$$\{w^\ast_1v_{2q_1+1},w^\ast_2v_{2q_1+2}\}, 
\{w^\ast_1v_{2q_1+3},w^\ast_2v_{2q_1+4}\}, \ldots, \{w^\ast_1v_{2q-1},w^\ast_2v_{2q}\} 
\subset E(G) \setminus E(B),$$ 
where $\{w^\ast_1,w^\ast_2\} \subseteq W$ and 
$u_i \in V(G_1)$, $v_i \in V(G_2)$ for $2q_1+1 \leq i \leq 2q$.

\bigskip

\noindent {\bf Proof of Claim 4:}
Suppose that $q_2 > 0$, i.e., $q_1 < q$.
Then we have $q_1 = \max\{0,|M_H|-|M_B|-\max\{p,q\}\}$ which 
means that $|M_H| \leq q_1+|M_B|+\max\{p,q\}$.
Thus,  
\begin{equation}
\max\{|U_1|,|U_2|\} \leq q_1+k-k''+\max\{p,q\}. \label{Claim4-0}
\end{equation}
Assume that $k' \leq 1$.
Since $\delta(G) \geq \frac{n}{2}$, 
for every vertex $v \in V(G_1)$, $N_G(v) \cap V(G_2) \neq \emptyset$, i.e., $U_1 = V(G_1)$.
From (\ref{Claim4-0}) and the lower bound in (\ref{G1}), 
it follows that 
$$\frac{n}{2}-k'-p+1 \leq  q_1+k-k''+\max\{p,q\},$$
that is,
$$n \leq 2(k+k'-k'')+2p+2q_1+2\max\{p,q\}-2.$$
Since $k'' \geq \frac{k'}{2}$, $2(k'-k'') \leq k' \leq k-1$.
Thus, we obtain 
$$n \leq 3k+2p+2q+2\max\{p,q\}-5,$$
which contradicts (\ref{lower-n1}).
Therefore, we have $k' \geq 2$.

Let 
$$W' = \left\{ 
w \in W\ |\ \min\{|N_G(w) \cap V(G_1)|, |N_G(w) \cap V(G_2)|\} \leq q_2 
\right\}.$$
Assume that $|W'| \geq k'-1$.
Let $W'_1 = W' \cap W_1$ and $W'_2 = W' \cap W_2$.
Then, for any $w \in W'_1$ (respectively, $w \in W'_2$), 
it holds that 
$|N_G(w) \cap V(G_2)| \leq q_2$ (respectively, $|N_G(w) \cap V(G_1)| \leq q_2$).

\medskip

\noindent Case 1: $|V(G_1) \cup W'_1| = |V(G_2) \cup W'_2|$.

Since $\delta(G) \geq \frac{n}{2}$ and 
$|W \setminus W'| \leq 1$, for any vertex $v$ in $V(G_1) \cup W'_1$,
$N_G(v) \cap (V(G_2) \cup W'_2) \neq \emptyset$.
This means that 
for any $v \in V(G_1) \setminus U_1$, $N_G(v) \cap W'_2 \neq \emptyset$. 
Therefore, we have $|V(G_1) \setminus U_1| \leq q_2|W'_2|$.
Since $|V(G_1)| \leq |V(G_2)|$, $|W'_1| \geq |W'_2|$. 
Thus it holds that $|W'_2| \leq \frac{|W'|}{2} \leq \frac{k'}{2}$.
Therefore, 
\begin{equation}
|V(G_1)| \leq |U_1|+\frac{q_2k'}{2}. \label{Claim4'-1}
\end{equation}
It follows from (\ref{G1}) and (\ref{Claim4-0}) that 
$$\frac{n}{2}-k'-p+1 \leq q_1+k-k''+\max\{p,q\}+\frac{q_2k'}{2},$$
that is, 
$$n \leq 2(k+k'-k'')+2p+2q_1+2\max\{p,q\}-2+q_2k'.$$
Since $2q_1+q_2k' \leq kq - (k-2)q_1-q_2$, we have 
\begin{equation}
n \leq kq+3k+2p+2\max\{p,q\}-3-(k-2)q_1-q_2. \label{Claim4-1}
\end{equation}

\medskip

\noindent Case 2: $|V(G_1) \cup W'_1| < |V(G_2) \cup W'_2|$.

If $|V(G_1) \cup W'_1|+2 \leq |V(G_2) \cup W'_2|$, or 
$|V(G_1) \cup W'_1|+1 = |V(G_2) \cup W'_2|$ and $|W'| = k'$, then 
for any $v \in V(G_1) \setminus U_1$, $|N_G(v) \cap W'_2| \geq 2$. 
Thus, in such a case, we have $2|V(G_1) \setminus U_1| \leq q_2|W'_2| \leq q_2k'$, i.e., (\ref{Claim4'-1}) 
and similarly to Case 1, (\ref{Claim4-1}) is obtained.
If $|V(G_1) \cup W'_1|+1 = |V(G_2) \cup W'_2|$ and $|W'| = k'-1$, then 
$|V(G_1) \setminus U_1| \leq q_2|W'_2|$ and $|W'_2| \leq \frac{|W'|+1}{2}
\leq \frac{k'}{2}$, from which we again have (\ref{Claim4'-1}) and (\ref{Claim4-1}).

\medskip

\noindent Case 3: $|V(G_1) \cup W'_1| > |V(G_2) \cup W'_2|$.

Similarly to Case 2, 
if  $|V(G_1) \cup W'_1| \geq |V(G_2) \cup W'_2|+2$, or 
$|V(G_1) \cup W'_1| = |V(G_2) \cup W'_2|+1$ and $|W'| = k'$, then
we have $2|V(G_2) \setminus U_2| \leq q_2|W'_1|$. 
In such a case, 
\begin{equation}
2|V(G_2)| \leq 2|U_2|+q_2k', \label{Claim4'-2}
\end{equation}
and it follows from (\ref{G2}) and (\ref{Claim4-0}) that 
$$n-k' \leq 2(q_1+k-k''+\max\{p,q\})+q_2k',$$
that is, 
$$n \leq 2k+k'-2k''+2q_1+2\max\{p,q\}+q_2k',$$
which implies 
\begin{equation}
n \leq kq+3k+2\max\{p,q\}-1-k'-(k-2)q_1-q_2. \label{Claim4-2}
\end{equation}

If $|V(G_1) \cup W'_1| = |V(G_2) \cup W'_2|+1$ and $|W'| = k'-1$, 
then $|V(G_2) \setminus U_2| \leq q_2|W'_1|$.
Thus, in such a case, 
\begin{equation}
|V(G_2)| \leq |U_2|+q_2(k'-1), \label{Claim4'-3}
\end{equation}
and  it also follows from (\ref{G2}) and (\ref{Claim4-0}) that 
$$\frac{n-k'}{2} \leq q_1+k-k''+\max\{p,q\}+q_2(k'-1),$$
that is, 
$$n \leq 2k+k'-2k''+2q_1+2\max\{p,q\}+2q_2(k'-1).$$
Since $2q_1+2q_2(k'-1) \leq 2kq -2q-2q_1(k-2)-2q_2$, we have 
\begin{equation}
n \leq kq+(k-2)q+3k+2\max\{p,q\}-1-k'-2(k-2)q_1-2q_2.  \label{Claim4-3}
\end{equation}

Now suppose that $q_2 \geq 2$.
Then all the inequalities (\ref{Claim4-1}), (\ref{Claim4-2}) and (\ref{Claim4-3})
contradict (\ref{lower-n2}).
Note that for (\ref{Claim4-3}), we also use the fact $k' \geq 2$. 
Therefore, it follows that $|W'| \leq k'-2$.
Hence, there exist two vertices $w^\ast_1,w^\ast_2 \in W$
such that $|N_G(w^\ast_j) \cap V(G_i)| \geq q_2+1$ for each $i,j \in \{1,2\}$. 
Let 
$$\left\{ \begin{array}{l}
E_1 = \{ uv \in E(G)\ |\ u \in V(G_1), v \in \{w^\ast_1,w^\ast_2\}\}
\setminus E(B),\\[2mm]
E_2 = \{ uv \in E(G)\ |\ u \in V(G_2), v \in \{w^\ast_1,w^\ast_2\}\}
\setminus E(B).
\end{array} \right.$$
Note that each edge-induced subgraph $\langle E_i \rangle_G$ is a bipartite graph. 
By the definition of $B$, it holds that $|N_B(w^\ast_j) \cap V(G_i)| \leq 1$
for each $i,j \in \{1,2\}$. 
Thus, in each $\langle E_i \rangle_G$, 
$w^\ast_1$ and $w^\ast_2$ have degree at least $q_2$
and any other vertex has degree at most two.
Therefore, $|E(\langle E_i \rangle_G)| \geq 2q_2$ 
and $|E(\langle E_i \rangle_G)| - \Delta(\langle E_i \rangle_G) \geq q_2$
for each $i \in \{1,2\}$.
Hence, by Theorem \ref{max-DEP}, each $\langle E_i \rangle_G$ has $q_2$ edge-disjoint edge-pairs.

Suppose that $q_2 = 1$.
We then modify the definition of $W'$ by replacing $q_2$ with $q_2+1$, i.e.,
$$W' = \left\{ 
w \in W\ |\ \min\{|N_G(w) \cap V(G_1)|, |N_G(w) \cap V(G_2)|\} \leq q_2+1\right\}.$$
By this modification, the inequalities (\ref{Claim4-1}), (\ref{Claim4-2}) and (\ref{Claim4-3})
are weakened with additional terms $k'$, $k'$ and $2(k'-1)$, respectively, in the right-hand sides.
However, from the existence of terms $-(k-2)q_1$, $-k'$ and 
$-2(k-2)q_1$ in (\ref{Claim4-1}), (\ref{Claim4-2}) and (\ref{Claim4-3}), respectively,
we can see that the weakened inequality for (\ref{Claim4-2}) still contradicts (\ref{lower-n2}) 
and if $q_1 \geq 2$ (respectively, $q_1 \geq 1$), then 
the weakened inequality for  (\ref{Claim4-1}) (respectively, (\ref{Claim4-3})) 
also contradicts (\ref{lower-n2}). 
Note that $k \ge k'+1 \geq 3$. 
When $q_1 = 1$, the weakened 
inequality for (\ref{Claim4-1}) is as follows: 
$$n \leq 4k+2p+2\max\{p,2\}-2+k'.$$
When $q_1 = 0$, the weakened 
inequalities for (\ref{Claim4-1}) and (\ref{Claim4-3})
are also as follows:
$$\left\{ \begin{array}{l}
n \leq 4k+4p-4+k',\\ 
n \leq 5k+2p-7+k'.
\end{array} \right.$$
These inequalities now contradict (\ref{lower-n1}).
Therefore, for each $i,j \in \{1,2\}$, 
${\rm deg}_{\langle E_i \rangle_G}(w^\ast_j) \geq 2$, which means that
each $\langle E_i \rangle_G$ has an edge-pair; actually, it has two edge-disjoint edge-pairs.
$\Box$

\bigskip

Let $$Q_1 = \cup_{1 \leq i \leq q_1}\{u_{2i-1}v_{2i-1},u_{2i}v_{2i}\}$$
and
$$Q_2 = \cup_{1 \leq i \leq q_2}\{u_{2q_1+2i-1}w^\ast_1, u_{2q_1+2i}w^\ast_2,
w^\ast_1v_{2q_1+2i-1}, w^\ast_2v_{2q_1+2i}\}.$$

Note that $Q_1$ consists of $q_1$ edge-disjoint edge-pairs, while 
$Q_2$ consists of $q_2$ edge-disjoint pairs of vertex-disjoint paths of order 3, i.e., 
$\{ (u_{2q_1+2i-1}, w^\ast_1, v_{2q_1+2i-1}), (u_{2q_1+2i}, w^\ast_2, v_{2q_1+2i}) \}$ for 
$1 \leq i \leq q_2$.

\bigskip

\noindent {\bf Claim 5.}\   
There are $q$ edge-disjoint Hamiltonian cycles $H_1,H_2,\ldots,H_q$
in $G-E(B)$. 

\bigskip

\noindent {\bf Proof of Claim 5:}
We first show the existence of $q_1$ edge-disjoint Hamiltonian cycles $H_1,H_2,\ldots,H_{q_1}$
in $G-E(B) \cup Q_2$. 
Using the $q_1$ edge-disjoint edge-pairs in $Q_1$, 
it is sufficient to show the existence of $q_1$ edge-disjoint Hamiltonian paths 
$H_{i,1},H_{i,2},\ldots,H_{i,q_1}$ in $\langle V(G_i) \cup W_i \rangle_{G-E(B) \cup Q_2}$
for each $i \in \{1,2\}$ such that $H_{1,j}$ (respectively, $H_{2,j}$) 
connects $u_{2j-1}$ and $u_{2j}$ (respectively, $v_{2j-1}$ and $v_{2j}$) 
for each $1 \leq j \leq q_1$.

Any $v \in V(G_1)$ (respectively, $V(G_2)$) has at least 
$\delta(G)-p-|W_2|$ (respectively, $\delta(G)-p-|W_1|$) neighbors in 
$\langle V(G_1) \cup W_1 \rangle_G$ (respectively, $\langle V(G_2) \cup W_2 \rangle_G$).
Moreover, any $v \in V(G_1)$ (respectively, $V(G_2))$ is incident with at most one edge of $B$
and at most $\min\{q_2,2\}$ edges of $Q_2$. 
Thus, every vertex in $V(G_1)$ (respectively, $V(G_2)$) 
has at least $\delta(G)-1-p-\min\{q_2,2\}-|W_2|$ 
(respectively, $\delta(G)-1-p-\min\{q_2,2\}-|W_1|$) neighbors in 
$\langle V(G_1) \cup W_1 \rangle_{G-E(B) \cup Q_2}$
(respectively, $\langle V(G_2) \cup W_2 \rangle_{G-E(B) \cup Q_2}$).

By the definitions of $W_1$ and $W_2$, 
any vertex in $W_1$ (respectively, $W_2$) 
has at least $\left\lceil \frac{\delta(G)-|W_2|}{2} \right\rceil$ 
(respectively, $\left\lceil \frac{\delta(G)-|W_1|}{2} \right\rceil$) 
neighbors in $\langle V(G_1) \cup W_1 \rangle_G$ (respectively, $\langle V(G_2) \cup W_2 \rangle_G$).
If $w \in W_1 \cap V_{I}(P^\ast_i)$ (respectively, $W_2 \cap V_{I}(P^\ast_i)$) for some $P^\ast_i$
with $|V(P^\ast_i)| = 3$, 
then $w$ is incident with exactly one edge of $B$
in $\langle V(G_1) \cup W_1 \rangle_{G}$ (respectively, $\langle V(G_2) \cup W_2 \rangle_{G}$).
If $w,w' \in W \cap V_{I}(P^\ast_i)$ for some $P^\ast_i$ with $|V(P^\ast_i)| = 4$, 
then from the construction of $P^\ast_i$ in the proof of Claim 2, 
neither $\{w,w'\} \subseteq W_1$ nor $\{w,w'\} \subseteq W_2$; 
that is, $w_{i,1} \in W_1$ and $w_{i,2} \in W_2$. 
Thus, every vertex in $W_1$ (respectively, $W_2$) is incident with exactly one edge of $B$ in 
$\langle V(G_1) \cup W_1 \rangle_{G}$ (respectively, 
$\langle V(G_1) \cup W_2 \rangle_{G}$).
Moreover, any vertex $w$ in $W_1$ (respectively, $W_2$) is incident with either $q_2$ edges of $Q_2$
or no edge of $Q_2$
in $\langle V(G_1) \cup W_1 \rangle_G$ (respectively, $\langle V(G_2) \cup W_2 \rangle_G$)
depending on whether $w \in \{w^\ast_1,w^\ast_2\}$ or not.
Therefore, every vertex in $W_1$ (respectively, $W_2$) has at least 
$\left\lceil \frac{\delta(G)-|W_2|}{2} \right\rceil-1-q_2$ 
(respectively, $\left\lceil \frac{\delta(G)-|W_1|}{2} \right\rceil-1-q_2$)
neighbors in $\langle V(G_1) \cup W_1 \rangle_{G-E(B) \cup Q_2}$ 
(respectively, $\langle V(G_2) \cup W_2 \rangle_{G-E(B) \cup Q_2}$). 

Now applying Lemma \ref{H} to 
$\langle V(G_i) \cup W_i \rangle_{G-E(B) \cup Q_2}$ for each $i \in \{1,2\}$, 
if the following inequalities hold, then we have desired $q_1$ edge-disjoint Hamiltonian paths: 
$$\left\{\begin{array}{l}
\delta(G)-1-p-\min\{q_2,2\}-|W_2|\geq \left\lceil \frac{|V(G_1) \cup W_1|}{2} \right\rceil +2q_1-1, \\[2mm]
\frac{\delta(G)-|W_2|}{2}-1-q_2 \geq |W_1|+2q_1,
\end{array} \right.$$
$$\left\{ \begin{array}{l}
\delta(G)-1-p-\min\{q_2,2\}-|W_1|\geq \left\lceil \frac{|V(G_2) \cup W_2|}{2} \right\rceil +2q_1-1, \\[2mm]
\frac{\delta(G)-|W_1|}{2}-1-q_2 \geq |W_2|+2q_1, 
\end{array} \right.$$
that is,
$$\left\{\begin{array}{l}
\delta(G) \geq p+2q_1+\min\{q_2,2\}+|W|+ \max\left\{ \left\lceil \frac{|V(G_1)|-|W_1|}{2} \right\rceil, 
\left\lceil \frac{|V(G_2)|-|W_2|}{2} \right\rceil \right\}, \\[2mm]
\delta(G) \geq 2q+2q_1+2+|W|+\max\{|W_1|,|W_2|\}.
\end{array} \right.$$
Since $\left\lceil \frac{|V(G_2)|}{2} \right\rceil \leq 
\left\lceil \frac{\left\lfloor \frac{n}{2} \right\rfloor +p -1}{2} \right\rceil 
\leq \frac{n}{4}+\frac{p}{2}$, the following is 
sufficient for the above inequalities: 
$$\frac{n}{2} \geq \max\left\{\frac{3p}{2}+2q_1+\min\{q_2,2\}+k'+\frac{n}{4},\ 2q+2q_1+2+2k'\right\},$$
that is,
$$n \geq 4k'+\max\{ 6p+8q_1+4\min\{q_2,2\}, 4q+4q_1+4\},$$
which follows from (\ref{lower-n1}). 

We next show the existence of $q_2$ edge-disjoint Hamiltonian cycles $H_{q_1+1},H_{q_1+2},\ldots,H_{q}$
in $G-\cup_{1 \leq i \leq q_1}E(H_i) \cup E(B)$. 
Let $W'' = W \setminus \{w^\ast_1,w^\ast_2\}$, 
$$W''_1 = \left\{ w \in W''\ \left|\ |N_{G-\cup_{1 \leq i \leq q_1}E(H_i) \cup E(B)}(w) \cap V(G_1)| 
\geq |N_{G-\cup_{1 \leq i \leq q_1}E(H_i) \cup E(B)}(w) \cap V(G_2)| \right. \right\}$$
and $W''_2 = W'' \setminus W''_1$.
Using the $q_2$ edge-disjoint pairs of vertex-disjoint paths of order 3 in $Q_2$, 
it is sufficient to show the existence of $q_2$ edge-disjoint Hamiltonian paths $H_{i,q_1+1},H_{i,q_1+2}$, 
$\ldots,H_{i,q}$
in $\langle V(G_i) \cup W''_i \rangle_{G-\cup_{1 \leq i \leq q_1}E(H_i) \cup E(B)}$ for each $i \in \{1,2\}$
such that $H_{1,j}$ (respectively, $H_{2,j}$) connects $u_{2j-1}$ and $u_{2j}$
(respectively, $v_{2j-1}$ and $v_{2j}$) for each $q_1+1 \leq j \leq q$.

Every $v \in V(G_1)$ (respectively, $V(G_2)$) has at least 
$\delta(G)-1-p-2q_1-|W''_2 \cup \{w^\ast_1,w^\ast_2\}|$ 
(respectively, $\delta(G)-1-p-2q_1-|W''_1 \cup \{w^\ast_1,w^\ast_2\}|$) neighbors in 
$\langle V(G_1) \cup W''_1 \rangle_{G-\cup_{1 \leq i \leq q_1}E(H_i) \cup E(B)}$ 
(respectively, $\langle V(G_2) \cup W''_2 \rangle_{G-\cup_{1 \leq i \leq q_1}E(H_i) \cup E(B)}$).
On the other hand, 
every vertex in $W''_1$ (respectively, $W''_2$) has at least 
$\left\lceil \frac{\delta(G)-2-2q_1-|W''_2\cup \{w^\ast_1,w^\ast_2\}|}{2} \right\rceil$
(respectively, $\left\lceil \frac{\delta(G)-2-2q_1-|W''_1\cup \{w^\ast_1,w^\ast_2\}|}{2} \right\rceil$)
neighbors in 
$\langle V(G_1) \cup W''_1 \rangle_{G-\cup_{1 \leq i \leq q_1}E(H_i) \cup E(B)}$ 
(respectively, $\langle V(G_2) \cup W''_2 \rangle_{G-\cup_{1 \leq i \leq q_1}E(H_i) \cup E(B)}$).

Applying Lemma \ref{H} to $\langle V(G_i) \cup W''_i \rangle_{G-\cup_{1 \leq i \leq q_1}E(H_i) \cup E(B)}$
for each $i \in \{1,2\}$, 
if the following inequalities hold, then we have desired $q_2$ edge-disjoint Hamiltonian paths: 
$$\left\{\begin{array}{l}
\delta(G)-1-p-2q_1-|W''_2 \cup \{w^\ast_1,w^\ast_2\}| \geq 
\left\lceil \frac{|V(G_1) \cup W''_1|}{2} \right\rceil +2q_2-1, \\[2mm]
\frac{\delta(G)-2-2q_1-|W''_2\cup \{w^\ast_1,w^\ast_2\}|}{2} \geq |W''_1|+2q_2,
\end{array} \right.$$
$$\left\{ \begin{array}{l}
\delta(G)-1-p-2q_1-|W''_1\cup \{w^\ast_1,w^\ast_2\}| \geq 
\left\lceil \frac{|V(G_2) \cup W''_2|}{2} \right\rceil +2q_2-1, \\[2mm]
\frac{\delta(G)-2-2q_1-|W''_1\cup \{w^\ast_1,w^\ast_2\}|}{2} \geq |W''_2|+2q_2,
\end{array} \right.$$
that is,
$$\left\{\begin{array}{l}
\delta(G) \geq p+2q+|W|+ \max\left\{ 
\left\lceil \frac{|V(G_1)|-|W''_1|}{2} \right\rceil,
\left\lceil \frac{|V(G_2)|-|W''_2|}{2} \right\rceil \right\},\\[2mm]
\delta(G) \geq 2q+2q_2+2+|W|+\max\{|W''_1|,|W''_2|\}.
\end{array} \right.$$
The following is sufficient for the above inequalities: 
$$\frac{n}{2} \geq \max\left\{ \frac{3p}{2}+2q+k'+\frac{n}{4},\ 2q+2q_2+2k'\right\},$$
that is,
\begin{equation}
n \geq 4k'+6p+8q. \label{Claim5'-1}
\end{equation}
This inequality follows from (\ref{lower-n1}). 
Hence, Claim 5 holds.
$\Box$

\bigskip

Let $G'' = G-\cup_{1 \leq i \leq q}E(H_i)$, 
$G''_1 = \langle V(G_1) \rangle_{G''}$ and 
$G''_2 = \langle V(G_2) \rangle_{G''}$.
By showing the following claim, the proof of Theorem \ref{main0} is complete.

\bigskip

\noindent {\bf Claim 6.}\  
$G''$ is $k$-connected.

\bigskip

\noindent {\bf Proof of Claim 6:}
If the following holds where $i \in \{1,2\}$, then by Theorem \ref{CH}, $G''_i$ is $k$-connected: 
\begin{equation}
\delta(G''_i) \geq \frac{|V(G''_i)|+k-2}{2}.  \label{C4}
\end{equation}
Since for each $i \in \{1,2\}$, $\frac{n}{2}+p-1 \geq |V(G''_i)|$ 
and $\delta(G''_i) \geq \delta(G)-p-2q-|W|$, 
the following is sufficient for (\ref{C4}):
$$\frac{n}{2}-p-2q-k' \geq \frac{\frac{n}{2}+p-1+k-2}{2},$$
that is, 
\begin{equation}
n \geq 2k+4k'+6p+8q-6, \label{Claim6'-1}
\end{equation}
which follows from (\ref{lower-n1}). 
Therefore, $G''_1$ and $G''_2$ are both $k$-connected.

By Claim 5, $B \subset G''$.
If the following holds, then every vertex $w \in W$ in $G''$ has at least $k$ neighbors in 
$V(G''_1) \cup V(G''_2)$: 
$$\delta(G'') - (|W|-1) \geq k.$$
This inequality holds if
$$\frac{n}{2}-2q-k'+1 \geq k,$$
that is,
\begin{equation}
n \geq 2k+2k'+4q-2, \label{Claim6'-2}
\end{equation}
which is satisfied by (\ref{lower-n1}).
Hence, it follows from Lemma \ref{k-connect} that  
$G''$ is $k$-connected.
$\Box$
\hfill $\blacksquare$

\bigskip

By setting $p = 2$ and $q = 1$ in Theorem \ref{main0}, we have Theorem \ref{main1}.
Note that we may let $H$ be a Hamiltonian cycle in $G$. 
Moreover, by setting $p = 2$ and $q = \ell \geq 2$ in Theorem \ref{main0}, 
we have the following corollary.

\begin{corollary} \label{main-cor1}
Let $G$ be a $k$-connected graph of order 
$n \geq \max\{k\ell+\max\{k\ell-6, 2\ell\}+3k,\ 6k+8\ell+4-2\left\lceil \frac{k-2}{k} \right\rceil \}$ 
with $\delta(G) \geq \frac{n}{2}$, where $k \geq 2$ and $\ell \geq 2$. 
If $G$ has a Hamiltonian cycle $H$ such that $G-E(H)$ is not $k$-connected, 
then there are $\ell$ edge-disjoint Hamiltonian cycles $H_1,H_2,\ldots,H_\ell$ in $G$ such that
$G-\cup_{1\leq i \leq \ell}E(H_i)$ is $k$-connected.
\end{corollary}

Furthermore, by setting $p = 2q = 2\ell$ for $\ell \geq 2$ in Theorem \ref{main0}, 
the following is obtained. 

\begin{corollary} \label{main-cor2}
Let $G$ be a $k$-connected graph of order 
$n \geq \max\{k\ell+\max\{k\ell,6\ell+2 \}+3k+2\ell-6,\ 6k+20\ell-8-2\left\lceil \frac{k-2}{k} \right\rceil \}$ 
with $\delta(G) \geq \frac{n}{2}$, where $k \geq 2$ and $\ell \geq 2$. 
If $G$ has $\ell$ edge-disjoint Hamiltonian cycles,
then there are $\ell$ edge-disjoint Hamiltonian cycles $H_1,H_2,\ldots,H_\ell$ in $G$ such that
$G-\cup_{1\leq i \leq \ell}E(H_i)$ is $k$-connected.
\end{corollary}

Applying Theorem \ref{N} to Corollary \ref{main-cor2}, 
Theorem \ref{main2} is obtained.

\section{Proof of Theorem \ref{main3}}

Let $G$ be a graph of order $n$ with $\kappa(G) = k \geq 2$ and $\delta(G) \geq \frac{n}{2}$. 
Let $W = \{w_1,w_2,\ldots,w_{k}\} \subset V(G)$
such that $G-W$ is disconnected.
In what follows, we employ similar notations and definitions used in Section 4
unless otherwise stated. 
Without introducing the spanning subgraph $H$,
we simply consider the existence of desired edge-disjoint Hamiltonian cycles in $G$.
That is, by modifying the proof in Section 4 with the setting
that $G' = G$, $k' = k$, $p = 0$ and $q = \ell$,  
we present lower bounds on $n$ for the existence of
$\ell$ edge-disjoint Hamiltonian cycles $H_1,H_2,\ldots,H_\ell$ such that
$\kappa(G-\cup_{1 \leq i \leq \ell}E(H_i)) = k$.

The following follows from (\ref{Claim1'-1}):
$$n \leq 4k-6.$$
Thus, we have Claim 1'.

\bigskip

\noindent {\bf Claim 1'.}\ 
If $n \geq 4k-5$, then 
$G-W$ has exactly two connected components $G_1$ and $G_2$
with $|V(G_1)| \leq |V(G_2)|$ such that 
$\frac{n}{2}-k+1 \leq |V(G_1)| \leq \frac{n-k}{2} \leq |V(G_2)| \leq \frac{n}{2}-1$. 

\bigskip

From a similar argument in the first paragraph of the proof of Claim 2, 
the following claim is obtained.

\bigskip

\noindent {\bf Claim 2'.}\ 
There are $k$ vertex-disjoint paths $P^\ast_1, P^\ast_2, \ldots , P^\ast_k$
of order 3 such that $P^\ast_i = (x_i,w_i,y_i)$ for $1 \leq i \leq k$, where 
$\{x_1,x_2,\ldots,x_{k}\} \subset V(G_1)$ and 
$\{y_1,y_2,\ldots,y_{k}\} \subset V(G_2)$.

\bigskip

From Claim 2', we have $k'' = k$ and $M_B = \emptyset$.
Note that $U_1 = U_2 = \emptyset$.
Since $H$ is not defined, we regard $M_H$ as the empty set. 
Thus, $q_1$ and $q_2$ are defined to be $0$ and $\ell$, respectively. 
Since $q_1 = 0$, we do not have to consider the claim corresponding to Claim 3.
The argument in the first paragraph of the proof of Claim 4 is also not needed
since $k' = k \geq 2$. 
The following inequalities are obtained 
from (\ref{Claim4'-1}), (\ref{Claim4'-2}) and (\ref{Claim4'-3}) with the lower bounds
in Claim 1': 
$$\left\{ \begin{array}{l}
n \leq k\ell+2k-2,\\ 
n \leq k\ell+k,\\
n \leq 2k\ell+k-2\ell. 
\end{array} \right.$$
In the case that $\ell = 1$ with the modified definition of $W'$, 
the inequalities (\ref{Claim4'-1}), (\ref{Claim4'-2}) and (\ref{Claim4'-3}) with 
the additional terms $k$, $k$ and $2(k-1)$, respectively, in the right-hand sides are 
expressed as follows:
$$\left\{ \begin{array}{l}
n \leq 4k-2,\\ 
n \leq 3k,\\
n \leq 5k-4. 
\end{array} \right.$$
From these inequalities, we have Claim 4'.

\bigskip 

\noindent {\bf Claim 4'.}\ 
If $n \geq \max\{2k\ell+ k-2\ell+1,5k-3\}$, then there are 
$\ell$ edge-disjoint edge-pairs 
$\{u_{2i-1}w^\ast_1,u_{2i}w^\ast_2\} \subset E(G) \setminus E(B)$ 
and $\ell$ edge-disjoint edge-pairs 
$\{w^\ast_1v_{2i-1},w^\ast_2v_{2i}\} \subset E(G) \setminus E(B)$
for $1 \leq i \leq \ell$, 
where $\{w^\ast_1,w^\ast_2\} \subseteq W$ and 
$u_i \in V(G_1)$, $v_i \in V(G_2)$ for $1 \leq i \leq 2\ell$.

\bigskip

The following follows from (\ref{Claim5'-1}):
$$n \geq 4k+8\ell.$$
Thus, we obtain Claim 5'.

\bigskip

\noindent {\bf Claim 5'.}\   
If $n \geq \max\{2k\ell+ k-2\ell+1,5k-3, 4k+8\ell\}$, 
then there are $\ell$ edge-disjoint Hamiltonian cycles $H_1,H_2,\ldots,H_\ell$ in $G-E(B)$. 

\bigskip

The following inequalities are obtained from (\ref{Claim6'-1}) and (\ref{Claim6'-2}): 
$$\left\{ \begin{array}{l}
n \geq 6k+8\ell-6,\\
n \geq 4k+4\ell-2. 
\end{array} \right.$$
Therefore, we have Claim 6'. 

\bigskip

\noindent {\bf Claim 6'.}\  
If $n \geq \max\{2k\ell+ k-2\ell+1, 4k+8\ell, 6k+8\ell-6\}$, then $G''$ is $k$-connected.

\bigskip

Consequently, if the following inequalities hold, then $G$ has desired Hamiltonian cycles: 
$$n \geq 2k\ell+k-2\ell+1$$
and
$$n \geq 6k+8\ell-4-2\left\lceil \frac{k-2}{k} \right\rceil.
$$
Hence, Theorem \ref{main3} holds.

\section{Concluding Remarks}

Based on our proofs, 
for any $k$-connected Dirac graph $G$ of order $n$ satisfying the lower bound on $n$ 
in Theorem \ref{main1} (respectively, \ref{main3}), 
a connectivity preserving Hamiltonian cycle
(respectively, connectivity preserving edge-disjoint Hamiltonian cycles) 
in $G$ can actually be constructed. 
This follows from the facts that there are algorithms for finding the structures 
whose existences are used in our proofs such as a minimum vertex-cut, vertex-disjoint paths,
a maximum set of edge-disjoint edge-pairs, and edge-disjoint Hamiltonian paths. 
In particular, an algorithm for finding a maximum set of edge-disjoint edge-pairs 
can be obtained from the proofs of Theorems \ref{decom} and \ref{max-DEP}. 
Furthermore, 
the $(n+1)$-closure $C_{n+1}(G)$ of a graph $G$ of order $n$ 
satisfying the condition in Theorem \ref{CKK} is indeed 
a complete graph, where $C_{n+1}(G)$ is obtained from $G$ by iteratively 
joining nonadjacent vertices whose degree sum is at least $n+1$,
and in such a case, a Hamiltonian path connecting any pair of vertices in $G$
can be found by applying the algorithm due to Bondy and Chv\'{a}tal \cite{BC}
for finding a Hamiltonian cycle in a graph of order $n$ whose $n$-closure is a complete graph. 

We have presented an $O(k)$ (respectively, $O(k\ell)$) lower bound on $n$
for the existence of a connectivity preserving Hamiltonian cycle (respectively,
connectivity preserving $\ell$ edge-disjoint Hamiltonian cycles) 
in a $k$-connected Dirac graph of order $n$.
While our $O(k)$ lower bound is asymptotically optimal in the sense that 
any $k$-connected graph has order $n = \Omega(k)$, 
it seems in general to be difficult to determine the tight lower bound on $n$
(which can be easily determined to be 7 when $\kappa(G) = 2$, however). 
It would be interesting to study 
whether our $O(k\ell)$ lower bound can be improved to $O(k+\ell)$. 
If we replace the lower bound of $\frac{n}{2}$ on $\delta(G)$ with $\frac{n+k-2}{2}$,
then such an improvement is possible.
Since every graph $G$ of order $n$ with $\delta(G) \geq \frac{n+k-2}{2}$ is $k$-connected, 
we can apply similar discussions in the proof of Theorem \ref{main0}.
Under the same definitions in Section 4,
it follows from $\delta(G) \geq \frac{n+k-2}{2}$ 
that any vertex $u \in V(G_1)$ has a neighbor $v \in V(G_2)$ in the spanning subgraph $H$,
i.e., $U_1 = V(G_1)$. 
This implies that $q_1 = q$, i.e., $q_2 = 0$.
Thus, we do not have to use Claim 4.
Note that the $O(k\ell)$ part of the lower bound on $n$ in Theorem \ref{main2} 
depends only on the lower bound condition
on $n$ required for Claim 4 in Theorem \ref{main0}.
Hence, the following holds.

\begin{theorem} \label{conc1}
For $k \geq 2$ and $\ell \geq 2$, 
there exists a funciton $f(k,\ell) = O(k+\ell)$ such that
every graph $G$ of order $n \geq f(k,\ell)$ with $\delta(G) \geq \frac{n+k-2}{2}$
contains $\ell$ edge-disjoint Hamiltonian cycles 
$H_1,H_2,\ldots,H_\ell$ such that 
$G-\cup_{1 \leq i \leq \ell}E(H_i)$ is $k$-connected.
\end{theorem}

\bigskip

\subsection*{Acknowledgments}

This work was supported by JSPS KAKENHI Grant Number JP19K11829.


\bigskip


\begin{thebibliography}{99}


\bibitem{BC} J.A. Bondy and V. Chv\'{a}tal, A method in graph theory, {\em Discrete Math.} 
\textbf{15} (1976) 111-136.

\bibitem{CH} G. Chartrand and F. Harary, Graphs with prescribed connectivities,
in: P. Erd\H{o}s and G. Katona (Eds.), Theory of Graphs, Academic Press, 1968, pp. 61--63.

\bibitem{CKK} G. Chartrand, S.F. Kapoor, and H.V. Kronk, A generalization of hamiltonian-connected graphs,
{\em J. Math. Pures Appl.} \textbf{48} (1969) 109--116. 

\bibitem{CL} G. Chartrand and L. Lesniak, Graphs \& Digraphs, 3rd edition, Chapman \& Hall, 1996. 

\bibitem{D} G.A. Dirac, Some theorems on abstract graphs,
{\em Proc. London Math. Soc.} \textbf{2} (1952) 69--81.

\bibitem{H} T. Hasunuma, Connectivity preserving trees in $k$-connected or $k$-edge-connected graphs,
{\em J. Graph Theory} \textbf{102} (2023) 423--435.

\bibitem{J} B. Jackson, Removal cycles in 2-connected graphs with minimum degree at least four,
{\em J. London Math. Soc} \textbf{21} (1980) 385--392.

\bibitem{M} W. Mader, Kreuzungsfreie $a$, $b$-Wege in endlichen Graphen, 
{\em Abhandlungen Math. Sem. Univ. Hamburg} \textbf{42} (1974) 187--204.

\bibitem{N} C.St.J.A. Nash-Williams, Edge-disjoint Hamiltonian circuits in graphs with vertices 
of large valency, 
in: L. Mirsky (Ed.), Studies in Pure Mathematics, Academic Press, 1971, pp. 157--183.

\bibitem{P} L. P\'{o}sa, A theorem concerning Hamilton lines, 
{\em Magyar. Tud. Akad. Mat. Kutat\'{o} Int. K\"{o}zl.} \textbf{7} (1962) 225--226.

\bibitem{R} S. Ruiz, Randomly decomposable graphs, 
{\em Discrete Math.} \textbf{57} (1985) 123--128.


\end{thebibliography}
\end{document}